\documentclass[12pt,reqno]{amsart}

\usepackage[latin1]{inputenc}
\usepackage[T1]{fontenc}
\usepackage[british]{babel}
\usepackage{eucal,url,relsize,xspace}

\theoremstyle{plain}
\newtheorem{theorem}{Theorem}[section]
\newtheorem{proposition}[theorem]{Proposition}
\newtheorem{lemma}[theorem]{Lemma}
\newtheorem{corollary}[theorem]{Corollary}

\theoremstyle{definition}
\newtheorem{definition}[theorem]{Definition}

\theoremstyle{remark}
\newtheorem{example}[theorem]{Example}
\newtheorem{remark}[theorem]{Remark}
\newtheorem*{acknowledgements}{Acknowledgements}

\numberwithin{equation}{section}

\DeclareMathOperator{\End}{End}
\DeclareMathOperator{\Id}{Id}
\DeclareMathOperator{\im}{Im}
\DeclareMathOperator{\sgn}{sgn}
\DeclareMathOperator{\HH}{\mathbb HH}
\DeclareMathOperator{\HP}{\mathbb HP}
\DeclareMathOperator{\RP}{\mathbb RP}
\newcommand{\hook}{\mathord{\mkern1mu\lrcorner\mkern1mu}}

\newcommand{\Lie}[1]{\operatorname{\textsl{#1}}}
\newcommand{\lie}[1]{\operatorname{\mathfrak{#1}}}
\newcommand{\Sl}{\lie{sl}}
\newcommand{\SP}{\Lie{Sp}}
\newcommand{\sP}{\lie{sp}}
\newcommand{\su}{\lie{su}}

\newcommand{\abs}[1]{\left\lvert #1\right\rvert}

\newcommand{\LC}{\nabla^{\text{LC}}}
\newcommand{\Nq}{\nabla^{\text{q}}}
\newcommand{\Ob}{\nabla^{\text{Ob}}}
\DeclareMathOperator{\UM}{\mathcal U}
\newcommand{\vol}{\text{vol}}

\newcommand{\HKT}{\textsmaller{HKT}\xspace}
\newcommand{\QKT}{\textsmaller{QKT}\xspace}

\newcommand{\NfB}{\( \mathcal N = 4\mathrm B \)\xspace}
\newcommand{\Dsym}[1][\alpha]{\( D(2,1;#1) \)\xspace}

\begin{document}
\title[Superconformal symmetry and HKT manifolds]{Superconformal symmetry
and hyperKähler manifolds with torsion}
\author{Yat Sun Poon}
\address[Poon]{Department of Mathematics\\
University of California at\linebreak Riverside\\
Riverside\\
CA 92521\\
USA}
\email{ypoon@math.ucr.edu}

\author{Andrew Swann}
\address[Swann]{Department of Mathematics and Computer Science\\
University of Southern Denmark\\
Campusvej 55\\
DK-5230 Odense M\\
Denmark}
\email{swann@imada.sdu.dk}

\begin{abstract}
  The geometry arising from Michelson \& Strominger's study of \( \mathcal
  N=4 \mathrm B \) supersymmetric quantum mechanics with superconformal \(
  D(2,1;\alpha) \)-symmetry is a hyperKähler manifold with torsion (HKT)
  together with a special homothety.  It is shown that different
  parameters~\( \alpha \) are related via changes in potentials for the HKT
  target spaces.  For \( \alpha\ne0,-1 \), we describe how each such HKT
  manifold \( M^{4m} \) is derived from a space~\( N^{4m-4} \) which is
  quaternionic Kähler with torsion and carries an Abelian instanton.
\end{abstract}

\subjclass{(2000) Primary 53C26; Secondary 53C07, 57S25, 81T60}
\keywords{HyperKähler, torsion, HKT, quaternionic Kähler, QKT, instanton,
symmetry}

\maketitle

\section{Introduction}

In the study of two-dimensional sigma models a variety of different
quaternionic geometries arise on the target spaces.  In the presence of a
Wess-Zumino term the metric connections have non-zero torsion.  For \NfB
rigid supersymmetry the target space carries an \HKT structure: the
geometry of a hyperKähler connection with totally skew symmetric torsion
\cite{Gibbons-PS:hkt-okt}.  For \NfB local symmetry the resulting
geometry~\cite{Howe-OP:QKT} is known as \QKT (quaternionic Kähler with
torsion).  The mathematical background of \HKT geometry was reported
in~\cite{Grantcharov-P:HKT}, where many examples were constructed.
Mathematical discussion of \QKT geometry may be found in \cite{Ivanov:QKT}.

Through the work of Maldacena~\cite{Maldacena:superconformal} there has
been much interest in field theories with superconformal symmetry.
Michelson \& Strominger \cite{Michelson-S:conformal} showed that for \NfB
rigid supersymmetry examples of quantum mechanical systems in one dimension
with actions of the superconformal groups~\Dsym may be obtained.  As
discussed in \cite{Britto-Pacumio-MSV:superconformal}, \Dsym has \(
\su(2)\oplus\su(2) \) as its algebra of R-symmetries and \Dsym[-2] is the
supergroup~\( \Lie{Osp}(4|2) \).  On the target space, Michelson \&
Strominger \cite{Michelson-S:conformal} show that the \HKT manifold
(locally) has a certain vector field~\( X \) generating one homothety and
three isometries, see equations~\eqref{eq:D21}. In this paper we
investigate the geometry of an \HKT manifold with such a vector field.  In
\cite{Poon-S:HKT-potential}, we showed that the length-squared of~\( X \)
gives a potential~\( \mu \) for the \HKT metric.  By transforming~\( \mu \)
we show in~\S\ref{sec:geography} that \Dsym-symmetries for different values
of~\( \alpha \) are closely related.  In particular, if an \HKT manifold
has a vector field~\( X \) generating a \Dsym-symmetry with \( \alpha<0 \)
and \( \alpha \ne -1 \), then the same manifold carries \HKT metrics with
\Dsym[\alpha']-symmetry for each \( \alpha'<0 \).  Similarly, any \(
\alpha'>0 \) may be obtained from any other~\( \alpha>0 \).

In~\S\ref{sec:quotient} we show that the vector fields generate an
infinitesimal action of the non-zero quaternions~\( \mathbb H^* \) and that
the quotient \( N^{4n}=M/\mathbb H^* \) carries a \QKT metric.  It turns
out, \S\ref{sec:base}, that this \QKT manifold comes equipped with an
instanton connection on its bundle \( \Lambda^{4n}T^*N \) of volume forms.
Locally \QKT metrics inducing instanton connections exist on any
quaternionic manifold, and from such a geometry in dimension~\( 4n \) we
construct in~\S\ref{sec:inverse} \HKT metrics with \Dsym-symmetry in
dimension~\( 4n+4 \).  As an interesting special case, we obtain \HKT
metrics with \Dsym[1]-symmetry over each quaternionic Kähler manifold of
negative scalar curvature.

Both the discussion of the parameter change for \Dsym-symmetry and the
bundle constructions relating \QKT and \HKT geometries naturally introduce
pseudo-Riemannian structures.  We therefore deal with \HKT geometry in this
generality from the outset.

If one sets the torsion to zero in this paper, then one recovers the
constructions of~\cite{Swann:MathAnn}, relating quaternionic Kähler
manifolds to hyperKähler manifolds with \Dsym[-2]-symmetry and hyperKähler
potentials.  This case is relevant to the discussion of superconformal
symmetry in \( \mathcal N=2 \) quantum
mechanics~\cite{DeWit-KV:hypermultiplets}.

\begin{acknowledgements}
  Andrew Swann is a member of \textsc{Edge}, Research Training Network
  \textsc{hprn-ct-\oldstylenums{2000}-\oldstylenums{00101}}, supported by
  The European Human Potential Programme.  He is grateful to the Department
  of Mathematics at the University of California at Riverside for
  hospitality during the initial stages of this work.  He wishes to thank
  Birte Feix and Stefan Ivanov for useful conversations, and Richard
  Cleyton for his comments on the manuscript.
\end{acknowledgements}

\section{Potentials and Superconformal Symmetry}
Let \( (M,g,I,J,K) \) be an \HKT manifold of dimension~\( 4m \) and
signature \( (4p,4q) \).  This means that \( I \), \( J \) and \( K \) are
integrable complex structures satisfying the quaternion identities, \( g
\)~is a hyper-Hermitian metric of signature~\( (4p,4q) \) and there is an
\( \SP(p,q) \)-connection~\( \nabla \) whose torsion tensor
\begin{equation*}
  c(X,Y,Z) = g(X,T(Y,Z))
\end{equation*}
is totally skew, where \( T(X,Y) = \nabla_XY - \nabla_YX - [X,Y] \).  The
integrability of~\( I \) implies
\begin{equation}
  \label{eq:T-I}
  T(IX,IY)-IT(IX,Y)-IT(X,IY)-T(X,Y)=0
\end{equation}
and that \( c \)~is of type \( (2,1)_I+(1,2)_I \).  Note that for a given
\( (g,I,J,K) \) there is at most one \HKT connection~\( \nabla \),
sometimes called the Bismut connection.

We set \( F_I(X,Y)=g(IX,Y) \) and define \( d_I \) on \( r \)-forms by
\begin{equation*}
  d_I\beta = (-1)^r IdI\beta,
\end{equation*}
where \( I\beta=\beta(-I\cdot,\dots,-I\cdot) \).  Similar forms and
operators are defined for \( J \) and~\( K \).  With these conventions the
torsion satisfies
\begin{equation*}
  -c=d_IF_I=d_JF_J=d_KF_K.
\end{equation*}

A \emph{potential} for an \HKT structure is a function such that
\begin{equation*}
  F_I = \tfrac12(dd_I+d_Jd_K)\mu, \quad\text{etc.}
\end{equation*}
Note that \( dd_I\mu = dId\mu \) and \( d_Jd_K\mu = -JdId\mu \).

In \cite[Corollary~4]{Grantcharov-P:HKT} it is shown that locally any
hypercomplex manifold~\( (M,I,J,K) \) admits a compatible \HKT metric with
potential.  On the other hand, Michelson \& Strominger
\cite[Appendix~C]{Michelson-S:conformal} show that for any open set of~\(
\mathbb R^{4m}=\mathbb H^m \) with the standard complex structures \( I \),
\( J \) and \( K \), any compatible \HKT metric admits a potential.  It is
an open question whether general \HKT structures admit potentials locally.
In \cite{Poon-S:HKT-potential} it was shown that those with \Dsym-symmetry
do.  Here we summarise that discussion.

Suppose we have an \HKT manifold with a vector field \( X \) satisfying
\begin{subequations}
  \label{eq:D21}
  \begin{align}
    \label{eq:Xg}
    &L_X g = a g\\
    \label{eq:IXg}
    &L_{IX} g = 0, \quad\text{etc.},\\
    \label{eq:IXI}
    &L_{IX} I = 0, \quad\text{etc.},\\
    \label{eq:IXJ}
    &L_{IX} J = bK, \quad\text{etc.},
  \end{align}
\end{subequations}
where \( a,b\in \mathbb R \) are constants and ``etc.''  means that the
versions of equations (\ref{eq:IXg}--\ref{eq:IXJ}) obtained by cyclically
permuting \( I \), \( J \) and \( K \) also hold.  By rescaling \( X \) we
may alter the constants \( a \) and~\( b \), but the point \( [a,b]\in
\RP(1) \cup \{(0,0)\} \) remains fixed.  We call a vector field~\( X \)
satisfying~\eqref{eq:D21} \emph{a special homothety of type~\( (a,b) \)}.

In the notation of \cite{Michelson-S:conformal}, \HKT geometry with this
symmetry arises from a quantum mechanical system with \Dsym-symmetry where
\( \alpha+1=a/b \).  For the standard flat metric on~\( \mathbb H^m=\mathbb
R^{4m} \), the vector field~\( X \) is given by dilation.  This case has \(
a=2 \), \( b=-2 \) and \( \alpha=-2 \).

For a vector field \( X \) satisfying equations \eqref{eq:Xg} and
\eqref{eq:IXg} on an \HKT manifold, it was shown
in~\cite{Poon-S:HKT-potential} that
\begin{equation*}
  \nabla X = \tfrac12 a\Id,
\end{equation*}
where \( \nabla \) is the torsion connection.  

\begin{lemma}
  \label{lem:exact}
  If \( a\ne0 \), then the one-form \( X^\flat \) is exact, \( g(X,X) \)~is
  non-constant and \( X\hook c=0 \).
\end{lemma}

\begin{proof}
  Since \( \nabla \) is a metric connection, we have
  \begin{equation*}
    d(g(X,X))= 2g(\nabla X,X) =  a g(\cdot,X) =  a X^\flat,
  \end{equation*}
  using the Lemma above.  This gives that \( X^\flat \) is exact and \(
  g(X,X) \)~is non-constant.  Now \( X^\flat \) is closed, so
  \begin{equation*}
    \begin{split}
      0 &= dX^\flat(Y,Z) = g(\nabla_YX,Z)-g(\nabla_ZX,Y)+g(X,T(Y,Z))\\
      &= c(X,Y,Z),
  \end{split}
\end{equation*}
  as required.
\end{proof}

Equations \eqref{eq:IXJ} and \eqref{eq:IXg} give
\begin{equation}
  \label{eq:IXc}
     IX\hook c - J(IX\hook c) = - (a+b) F_I
\end{equation}
which together with the Lemma is the basis for the proof of the following
result.

\begin{theorem}[\cite{Poon-S:HKT-potential}]
  \label{thm:mu}
  If \( X \) is a special homothety of type~\( (a,b) \) with \( a\ne 0,b
  \), then the function
  \begin{equation*}
    \mu = \frac2{a(a-b)} g(X,X) 
  \end{equation*}
  is an \HKT potential. \qed
\end{theorem}

\section{Parameter Changes}
\label{sec:geography}

Suppose \( M \)~is an \HKT manifold with potential~\( \mu \).  If \( f
\)~is a smooth function then Grantcharov \& Poon \cite{Grantcharov-P:HKT}
showed that
\begin{equation}
  \label{eq:gf}
  g_f = f'(\mu)\, g + \tfrac12f''(\mu)\, (d^{\mathbb H}\mu)^2
\end{equation}
is an \HKT metric with potential~\( f(\mu) \) whenever \( g_f \) is
non-degenerate, where \( (d^{\mathbb H}\mu)^2 = d\mu^2 + (Id\mu)^2 +
(Jd\mu)^2 + (Kd\mu)^2 \).  (Our conventions give slightly different
coefficients to those in~\cite{Grantcharov-P:HKT}.)  Such transformations
\( \mu \mapsto f(\mu) \) allow us to relate \HKT structures with different
\Dsym-symmetries, perhaps at the cost of changing the signature of the
metric.

\begin{proposition}
  Let \( (M,g,\nabla,X) \) be an \HKT manifold with a special homothety~\(
  X \) of type~\( (a,b) \), \( a\ne0,b \), and let \( \mu \)~be the
  potential found in Theorem~\ref{thm:mu}.  Suppose \( f \) is a smooth
  real-valued function.  Then \( f(\mu) \) is a potential for an \HKT
  metric~\( g_f \) with \( X \) a special homothety of type~\( (a',b') \)
  if either
  \begin{enumerate}
  \item[(a)] \( f(\mu)=\abs{\mu}^k \), \( k\ne 0,b/a \), \( (a',b')=(ka,b)
  \) or 
  \item[(b)] \( f(\mu)=\log\abs\mu \), \( (a',b')=(0,b) \).
  \end{enumerate}
  Up to homothety, these are the only possibilities for~\( g_f \).
\end{proposition}

Note that \( g_f \) is only defined away from the set \( \mu = 0 \).

\begin{proof}
  It is sufficient to determine when \( X \) is a homothety for the
  metric~\( g_f \) of~\eqref{eq:gf}.  As \( L_Xd\mu = a\,d\mu \), we have
  \begin{equation*}
    L_X g_f
    = a\left(\frac{f''}{f'}\mu+1\right) f'g +
    a\left(\frac{f'''}{f''}\mu+2\right) \tfrac12f''(d^{\mathbb H}\mu)^2.
  \end{equation*}
  So \( X \) is conformal only if 
  \begin{equation*}
    \left(\frac{f''}{f'}-\frac{f'''}{f''}\right)\mu = 1.
  \end{equation*}
  The left-hand side of this equation is \(
  -\mu\tfrac{d}{d\mu}\log\abs{f''/f'} \), so \( f''/f' = (k-1)\mu^{-1} \)
  for some~\( k \).  As \( f''/f'=\tfrac{d}{d\mu} \log\abs{f'} \),
  integration gives \( f' = A\abs{\mu}^{k-1} \).  The constant \( A \) only
  scales \( g_f \) by a constant, so we may take~\( A=\pm 1 \).  Note that
  \( X \) now scales~\( g_f \) by a constant.  Finally, we may integrate
  one more time to get the desired functions.

  From the form of~\( g_f \) we have
  \begin{gather*}
    g_{\abs\mu^k}(X,X) = \frac{ka-b}{a-b} k\sgn\mu\abs\mu^{k-1}g(X,X),\\
    g_{\abs\mu^k}(Y,Y) =k\sgn\mu\abs\mu^{k-1}g(Y,Y),
  \end{gather*}
  when \( Y \) is \( g \)-orthogonal to the quaternionic span of~\( X \).
  Thus \( ka-b \) needs to be non-zero to ensure a non-degenerate metric.
\end{proof}

When \( \mu \) is positive we see that \( g_{\mu^k} \) has the same
signature as \( g \) if and only if \( ka-b \) has the same sign as \( a-b
\).  For \( f(\mu)=\log\mu \) we have
\begin{equation*}
  g_{\log\mu}(X,X) = \frac b{(b-a)\mu}g(X,X),\qquad g_{\log\mu}(Y,Y) =
  \frac1\mu g(Y,Y),
\end{equation*}
which has the signature of~\( g \) only if \( b \) and \( b-a \) have the
same sign.  Recalling that \( \alpha = -1 + a/b \), we have:

\begin{corollary}
  Let \( M \) be a definite \HKT manifold with a nowhere zero special
  vector field generating a \Dsym-symmetry, with \( \alpha \)~finite.
  \begin{enumerate}
  \item \hfuzz=1.5pt If \( \alpha>0 \) then for each \( \alpha'>0 \) the
    hypercomplex manifold~\( M \) also admits a definite \HKT metric
    generating a \Dsym[\alpha']-symmetry.
  \item If \( \alpha<0 \) and \( \alpha\ne-1 \) then \( M \)~also admits a
    definite \HKT metric with \Dsym[\alpha']-symmetry for each \( \alpha'<0
    \).
  \end{enumerate}
\end{corollary}

\begin{proof}
  By rescaling we may choose the special homothety~\( X \) to have type~\(
  (a,1) \).  By replacing \( g \)~by \( -g \) if necessary, we may also
  ensure that \( \mu>0 \).  The transformations \( \mu\to \mu^k \) give the
  desired metrics for \( \alpha'\ne-1 \).  The case \( \alpha'=-1 \) is
  obtained as \( g_{\log\mu} \) when \( \alpha<0 \) and \( \alpha\ne-1 \).
\end{proof}

\begin{remark}
  Locally, one may change the parameter in a \Dsym[-1]-symmetry only if \(
  X\hook c = 0 \).  In this situation the distribution orthogonal to~\( X
  \) is integrable and we may locally solve the equation \( d\mu = \mu
  X^\flat \) to obtain an \HKT potential~\( \mu \).  This may be used to
  form \( g_{\mu^2} \) which has \( X \) as a special homothety of type \(
  (g(X,X),b) \).
\end{remark}

\section{The QKT Quotient}
\label{sec:quotient}
Let us first define what is meant by a \QKT structure on a manifold~\( N \)
of dimension~\( 4n \).  

The data consists of a metric~\( g \), a connection~\( \nabla^N \) and
subbundle~\( \mathcal G \) of \( \End TN \).  The bundle \( \mathcal G \)
should locally have a linear basis \( \{I_N, J_N, K_N\} \) satisfying the
quaternion relations \( I_N^2=J_N^2=K_N^2=-1 \) and \( I_NJ_N=K_N=-J_NI_N
\).  Call such a triple \( \{I_N, J_N, K_N\} \) a \emph{quaternion basis}
for~\( \mathcal G \).  The metric \( g \) is required to be Hermitian with
respect to each of these basis elements \( I_N \), \( J_N \) and \( K_N \).
The pair \( (g,\mathcal G) \) thus constitutes an almost quaternion
Hermitian manifold.

The connection \( \nabla^N \) should be metric, \( \nabla^Ng=0 \), and
quaternionic, so \( \nabla^NI_N \) is a linear combination of \( J_N \) and
\( K_N \).  In addition its torsion tensor \( c^N \) should be totally skew
and of type \( (2,1)+(1,2) \) with respect to each \( I_N \).  If these
conditions are satisfied then \( (g, \mathcal G, \nabla^N) \) is called a
\QKT structure.

Note that, the type condition is the same as saying that the torsion \( T^N
\)~satisfies the relation~\eqref{eq:T-I} for each choice of \( I=I_N \).

\begin{lemma}
  If \( X \) is a special homothety of type~\( (a,b) \) with \( a,b\ne0 \),
  then \( X, IX, JX, KX \) generate a local action of \( \mathbb H^* \).
\end{lemma}

\begin{proof}
  We compute the Lie brackets.  First,
  \begin{equation*}
    \begin{split}
      [X,IX] &= \nabla_X(IX) - \nabla_{IX}X - T(X,IX)\\
      &= \tfrac 12aIX - \tfrac12 aIX = 0,
    \end{split}
  \end{equation*}
  so \( X \) is central.  For the remainder we have
  \begin{equation*}
    \begin{split}
      [IX,JX] &= \nabla_{IX}JX- \nabla_{JX}IX - T(IX,JX)\\
      &= -aKX - T(IX,JX).
    \end{split}
  \end{equation*}
  However, equation~\eqref{eq:IXc} gives
  \begin{equation*}
    \begin{split}
      c(IX,JX,Z) &= -c(IX,X,JZ)-(a+b)F_I(JX,Z) \\
      &= -(a+b) (KX)^\flat(Z)
  \end{split}
\end{equation*}
  and hence \( [IX,JX]= bKX \).
\end{proof}

\begin{remark}
  When \( a = 0 \), one can show that \( X, IX, JX, KX \) generate a local
  action of~\( \mathbb H^* \) if \( X\hook c = 0 \).  In this situation one
  has \( \nabla X = 0 \), so \( [X,IX]=0 \) if and only if \( X\hook c=0
  \).  The relation \( [IX,JX]=-b(KX)^\flat \) follows from~\eqref{eq:IXc}.
\end{remark}

\begin{theorem}
  \label{thm:H-quotient}
  Let \( M^{4m} \) be an \HKT manifold with a special homothety \( X \) of
  type \( (a,b) \) with \( a,b,0 \) unequal and potential \( \mu \) from
  Theorem~\ref{thm:mu}.  Suppose that the vector fields \( IX,JX,KX \) are
  complete and let \( \SP(1) \) be the corresponding subgroup of \( \mathbb
  H^* \).  For \( x\ne0 \), the group \( \SP(1) \) acts semi-freely on \(
  \mu^{-1}(x) \) and the quotient \( \mu^{-1}(x)/\SP(1) \) is a \QKT
  orbifold.
\end{theorem}

The proof will occupy the rest of this section.  First note that as \( \mu=
2g(X,X)/(a(a-b)) \) and \( d\mu=2X^\flat/(a-b) \), each \( x\ne0 \) is a
regular value of~\( \mu \) and \( X \)~is not null on~\( \mu^{-1}(x) \).
The group \( \SP(1) \) acts semi-freely on
\begin{equation*}
  S=\mu^{-1}(x),
\end{equation*}
since \( IX \) preserves \( g \) and commutes with~\( X \).  The action of
\( \SP(1) \) is isometric by~\eqref{eq:IXg}, so we get a Riemannian metric
\( g_N \) on the quotient \( N=S/\SP(1) \).  Let \( \pi\colon S\to N \) be
the projection and write \( i\colon S\to M \) for the inclusion.

We define local almost complex structures \( I_N \), \( J_N \) and \( K_N
\) on \( N \) as follows.  Since \( \ker \pi_* \) is spanned by \( IX \),
\( JX \) and \( KX \), the horizontal distribution \( \mathcal H = (\ker
\pi_*)^\bot\subset TS \) is of dimension \( 4n \), where \( n=m-1 \) and is
preserved by \( I \), \( J \) and~\( K \).  Thus each point \( s\in
\pi^{-1}(p) \) defines a triple \( I_N,J_N,K_N \) of almost complex
structures on \( T_pN\cong \mathcal H_s \).  If \( s' \) is another point
of \( \pi^{-1}(p) \), then \( s'=gs \) for some \( g\in \SP(1) \).  But the
action of \( g \) permutes \( I \), \( J \) and \( K \), so the almost
complex structures~\( I_N',J_N',K_N' \) determined by~\( s' \) are linear
combinations of \( I_N \), \( J_N \) and \( K_N \).  The metric \( g_N \)
is Hermitian with respect to each of these almost complex structures.

In order to construct a \QKT structure on \( N \) we need a connection~\(
\nabla^N \).  On \( M \) we have \( \nabla = \LC+\tfrac12T \), where \( \LC
\) is the Levi-Civita connection of~\( g \).  This equation is also valid
on~\( S \), since Lemma~\ref{lem:exact} says that \( X\hook c=0 \), so \( T
\)~has no component normal to~\( S \).  In particular, \( S \) with the
induced metric naturally carries a metric connection whose torsion is skew.

Now on \( S \), the \( \SP(1) \)-action is isometric, and so preserves the
Levi-Civita connection.  On the other hand the torsion~\( T \) is \( \SP(1)
\)-invariant.  To see this, first note that for a two-form~\( \beta \) we
have
\begin{equation*}
  \begin{split}
    L_{IX}(J\beta)(Y,Z)
    &= (IX)(\beta(JY,JZ))-\beta(J[IX,Y],JZ)\\
    &\qquad-\beta(JY,J[IX,Z])\\
    &= (IX)(\beta(JY,JZ))-\beta([IX,JY]-bKY,JZ)\\
    &\qquad -\beta(JY,[IX,JZ]-bKZ)\\
    &= J(L_{IX}\beta)(Y,Z) + b\beta(KY,JZ) +b\beta(JY,KZ).
  \end{split}
\end{equation*}
If \( \beta \) is of type~\( (1,1)_I \) then this simplifies to \(
L_{IX}(J\beta) = JL_{IX}\beta \).  Taking \( \beta=dId\mu \), we now have
\begin{equation*}
  \begin{split}
    L_{IX}c
    &= -\tfrac12 L_{IX}d_Id_Jd_K\mu
    = \tfrac12 d_IL_{IX}JdId\mu \\
    &= -\tfrac12 d_Id_Jd_KL_{IX}\mu =0.
  \end{split}
\end{equation*}

We define \( \nabla^N \) by 
\begin{equation*}
  \nabla^N_AB = \pi_* \nabla_{\tilde A}\tilde B,
\end{equation*}
where \( \tilde A \) and \( \tilde B \) are the \( \SP(1) \)-invariant
lifts of \( A \) and \( B \) on~\( N \) to \( \mathcal H\subset TS \).  By
the above discussion, we have
\begin{equation*}
  \nabla^N_AB = \nabla^{\text{LC},N}_AB + \xi_AB,
\end{equation*}
where \( \xi_AB=\tfrac12\pi_*T(\tilde A,\tilde B) \) and \(
\nabla^{\text{LC},N} \) is the Levi-Civita connection of~\( g^N \).  Since
the torsion~\( T \) is \( \SP(1) \)-invariant, we see that \( \nabla^N \)
is well-defined.  Also the torsion three-tensor is
\begin{equation*}
  c^N(A,B,C)=c(\tilde A,\tilde B,\tilde C),
\end{equation*}
so \( c^N \) is skew symmetric and of type \( (2,1)+(1,2) \) for each
almost complex structure~\( I_N \).

Finally, we need to check that \( \nabla^N \) preserves the metric and
almost complex structures.  As \( I_NA = \pi_*((rI+sJ+tK)\tilde A) \) for
some functions \( r \), \( s \) and \( t \) on~\( S \), we have
\begin{equation*}
  \begin{split}
    (\nabla^N_A I_N)(B)
    & = \pi_*
    \left[
      \nabla_{\tilde A}(rI+sJ+tK)\tilde B - (rI+sJ+tK)\nabla_{\tilde
      A}\tilde B
    \right]\\
    & = \pi_* ((\tilde Ar)I + (\tilde As)J + (\tilde At)K )\tilde B,
  \end{split}
\end{equation*}
which is a linear combination of \( I_NB \), \( J_NB \) and \( K_NB \).
Thus \( \nabla^N \) is quaternionic.  Looking at~\( g^N \), we get
\begin{equation*}
  \begin{split}
    \pi^*(\nabla^N g^N)(A,B,C) &=
    \pi^*(Ag^N(B,C)-g^N(\nabla^N_AB,C)\\
    &\qquad-g^N(B,\nabla^N_AC)) \\
    &= (\nabla g)(\tilde A,\tilde B,\tilde C) = 0,
\end{split}
\end{equation*}
so \( \nabla^N \) is metric.

As the \( c^N \) has type \( (2,1)+(1,2) \) for each \( I_N \), we see that
\( N \) is \QKT. \qed

\section{The Geometry of the Quotient}
\label{sec:base}
Let \( N \) be the \QKT quotient constructed in the previous section.  Here
we will investigate special properties of this manifold.  First note that
if the \HKT metric has signature \( (4p,4q) \) then the metric on~\( N \)
has signature \( (4p-4,4q) \) if \( g(X,X)>0 \) or \( (4p,4q-4) \) if \(
g(X,X)<0 \).  In particular, the metric on~\( N \) may be definite even if
the original \HKT metric on~\( M \) is not.

Since \( \nabla^N \) is a quaternionic connection we have that the
curvature satisfies
\begin{equation}
  \label{eq:R-I}
  R^N_{A,B}I_N = - \beta_K(A,B)J_N + \beta_J(A,B)K_N,\quad\text{etc.}
\end{equation}
for some two-forms \( \beta_I \), \( \beta_J \) and \( \beta_K \).

\begin{proposition}
  \label{prop:sp1}
  \( \beta_I \) is of type~\( (1,1) \) with respect to \( I_N \).
\end{proposition}

\begin{proof}
  Write \( A=\pi_*\tilde A \) and \( I_NA = \pi_*((rI+sJ+tK)\tilde A) \) as
  the push-forward of invariant vector fields on~\( S \).  Then \( rI+sJ+tK
  \) is invariant under the action of \( IX \), \( JX \) and~\( KX \), and
  we have
  \begin{equation*}
    \begin{split}
      0 &= L_{IX}(rI+sJ+tK) \\
      &= ((IX)r) I + ((IX)s) J + ((IX)t) K + bsK - btJ,
    \end{split}
  \end{equation*}
  so \( (IX)r = 0 \), \( (IX)s = bt \) and \( (IX)t = - bs \).  Further
  such relations are obtained by considering the Lie derivative with
  respect to \( JX \) and \( KX \).
  
  The curvature \( R^N_{A,B}I_N \) pulls-back to
  \begin{equation*}
    \begin{split}
      ([\tilde A,\tilde B]^{\mathcal V}r)I &+ 
      ([\tilde A,\tilde B]^{\mathcal V}s)J + 
      ([\tilde A,\tilde B]^{\mathcal V}t)K
      \\
      &=\frac1{g(X,X)}
      \biggl\{
        \Bigl(
          g([\tilde A,\tilde B],IX)((IX)r) \\
      &\qquad\qquad\qquad    + g([\tilde A,\tilde B],JX)((JX)r) \\
      &\qquad\qquad\qquad    + g([\tilde A,\tilde B],KX)((KX)r)
        \Bigr)
        I + \text{etc.}
      \biggr\}
      \\
      &=\frac b{g(X,X)}
      \Bigl\{
        -g([\tilde A,\tilde B],JX)tI+g([\tilde A,\tilde B],KX)sI\\
        &\qquad\qquad\qquad
        +g([\tilde A,\tilde B],IX)tJ-g([\tilde A,\tilde B],KX)rJ\\ 
        &\qquad\qquad\qquad
        -g([\tilde A,\tilde B],IX)sK+g([\tilde A,\tilde B],JX)rK
      \Bigr\}.
    \end{split}
  \end{equation*}
  Evaluating at \( (r,s,t)=(1,0,0) \) we see that \( \beta_J \) is given by
  \begin{equation*}
    \begin{split}
      g([\tilde A,\tilde B],JX)
      &= g(\nabla_{\tilde A}\tilde B - \nabla_{\tilde B}\tilde A
      - T(\tilde A,\tilde B), JX)\\
      &= \tilde A g(\tilde B,JX) -g(\tilde B,\nabla_{\tilde A}JX) \\
      &\qquad - \tilde B g(\tilde A,JX) + g(\tilde A,\nabla_{\tilde B}JX)\\
      &\qquad - c(\tilde A,\tilde B,JX) \\
      &= -a F_J(\tilde A,\tilde B)-c(\tilde A,\tilde B,JX) \\
      &= \frac{b-a}2 dJd\mu (\tilde A,\tilde B),
    \end{split}
  \end{equation*}
  which is of type \( (1,1)_J \).
\end{proof}

To understand this curvature better we need to consider the relationship of
the torsion connection with the underlying quaternionic geometry.  

First recall that \HKT structures are built on top of hypercomplex
structures.  In \cite{Obata:connection} Obata showed that a hypercomplex
manifold admits a unique torsion-free connection~\( \Ob \) preserving the
complex structures.  Similarly, a \QKT manifold admits torsion-free
connections preserving the quaternionic structure; these are no longer
unique but form an affine space modelled on the one-forms.

\begin{lemma}
  \label{lem:xi}
  Suppose \( (Q,\mathcal G) \) is an almost quaternionic manifold and that
  \( T \in \Gamma(\Lambda^2 T^*Q\otimes TQ) \) is a tensor satisfying
  equation~\eqref{eq:T-I} with respect to almost complex structures \(
  I,J,K \) forming a quaternion basis of~\( \mathcal G \).  Define \( \xi
  \in \Gamma T^*Q\otimes \End(TQ) \) by
  \begin{equation*}
    \begin{split}
      \xi_YZ
      &= -\tfrac12 T(Y,Z)\\
      &\qquad + \tfrac16
      \bigl(
      IT(Y,IZ)+JT(Y,JZ)+KT(Y,KZ)\\
      &\qquad\qquad
      -IT(IY,Z)-JT(JY,Z)-KT(KY,Z)
      \bigr)\\
      &\qquad - \tfrac1{12}
      \bigl(
      IT(JY,KZ)+JT(KY,IZ)+KT(IY,JZ)\\
      &\qquad\qquad
      -IT(KY,JZ)-JT(IY,KZ)-KT(JY,IZ)
      \bigr).
    \end{split}
  \end{equation*}
  Then \( \xi \) satisfies
  \begin{enumerate}
  \item \( \xi_YZ-\xi_ZY=T(Y,Z) \),
  \item \( \xi I = 0 \),
  \item \( \xi \) is independent of the choice of quaternion basis \(
    \{I,J,K\} \) for~\( \mathcal G \).
    \qed
  \end{enumerate}
\end{lemma}

The proof of the Lemma and the following Proposition are straightforward
computations.

\begin{proposition}
  \label{prop:Nq}
  Let \( (N,\mathcal G,\nabla) \) be a \QKT manifold.  Set
  \begin{equation*}
    \Nq = \nabla + \xi ,
  \end{equation*}
  where \( \xi \)~is given by Lemma~\ref{lem:xi}.  Then \( \Nq \) is a
  torsion-free quaternionic connection on~\( N \) with \( \Nq A = \nabla
  A\), for each local section~\( A \) of~\( \mathcal G \).
  \qed
\end{proposition}

The uniqueness of the Obata connection now gives:

\begin{corollary}
  Let \( (M,\nabla) \) be an \HKT manifold.  Then the Obata connection
  is given by \( \Ob = \nabla + \xi \) with \( \xi \) defined in
  Lemma~\ref{lem:xi}.  
  \qed
\end{corollary}

The above Proposition is useful as it allows us to apply information about
quaternionic curvature from the work of Alekseevsky \& Marchiafava
\cite{Alekseevsky-M:subordinated} to give an interpretation of
Proposition~\ref{prop:sp1}.

First consider the case when \( \dim N > 4 \).  As \( \Nq \) is a
quaternionic connection, its curvature~\( R^{\,\text{q}} \) may be written
as
\begin{equation}
  \label{eq:curvature}
   R^{\,\text{q}} = R^B + W^{\text{q}}
\end{equation}
where \( W^{\text{q}} \) is an algebraic curvature tensor for \(
\Sl(n,\mathbb H) \) and \( R^B \) is determined by an element \( B \in
\Gamma(T^*N\otimes T^*N) \).  The component~\( W^{\text{q}} \) is
independent of the choice of the torsion-free quaternionic connection and
acts trivially both on~\( \mathcal G \) and the real canonical bundle~\(
\kappa^{\mathbb R}=\Lambda^{4n}T^*N \).  The curvature of~\( \mathcal G \)
is
\begin{equation*}
  \beta_I(Y,Z) = 2 \bigl( B(Y,I_NZ)-B(Z,I_NY) \bigr).
\end{equation*}
Write \( 2B=\lambda_I^{1,1}+\lambda_I^{2,0}+\sigma_I^{1,1}+\sigma_I^{2,0}
\) according to the splittings of \( \Lambda^2T^*N \) and \( S^2T^*N \)
with respect to~\( I \).  Here \( \lambda_I^{2,0} \) denotes the component
of the skew-symmetric part of~\( 2B \) in \(
\Lambda^{2,0}_IN+\Lambda^{0,2}_IN \), and similar notation is used for~\(
\sigma \).  We have \( \beta_I = \sigma_I^{1,1}(\cdot,I\cdot) +
\lambda^{2,0}_I(\cdot,I\cdot) \).  Proposition~\ref{prop:sp1} thus implies
that \( \lambda^{2,0}_I=0 \) for all~\( I \), i.e., the skew-part of~\( B
\) is of type \( (1,1) \) for all~\( I \).  However, the skew-part of~\( B
\) is the curvature of~\( \nabla^q \) on~\( \kappa^{\mathbb R} \).  Thus
the curvature of \( \kappa^{\mathbb R} \) is of type~\( (1,1) \) for all \(
I \).  In other words \( \Nq \) induces a instanton connection on~\(
\kappa^{\mathbb R} \).

For \( N \) a four-manifold, the decomposition~\eqref{eq:curvature} has an
extra term \( W_- \): the anti-self-dual part of the Weyl curvature.  If we
assume \( W_-=0 \), then the above analysis applies also in dimension four.

\begin{definition}
  Let \( N \) be a \QKT manifold.  If \( \dim N = 4 \) suppose in addition
  that \( N \) is self-dual.  We say that \( N \) is \emph{of instanton
  type} if \( \Nq \) induces an instanton connection on the real canonical
  bundle~\( \kappa^{\mathbb R} \).
\end{definition}

\begin{remark}
  Comparing with quaternionic or quaternionic Kähler geometry it would have
  been natural to include self-duality in the definition of \QKT manifolds
  for dimension four.  However, that goes against the established
  definitions in the \QKT literature.
\end{remark}

Note that the above discussion implies that \QKT manifolds of instanton
type are precisely those for which the curvature forms \( \beta_I \) are of
type~\( (1,1)_I \).  We summarise the above discussion in the following
result.

\begin{theorem}
  Let \( N \) be a \QKT manifold which is an \( \mathbb H^* \)-quotient of
  an \HKT manifold as in Theorem~\ref{thm:H-quotient}.  Then \( N \) is of
  instanton type.  \qed
\end{theorem}

This condition may related to the torsion of the \QKT manifold as follows.
Recall that there is a torsion \emph{one-form}~\( \tau \) given by
\begin{equation*}
  \tau(A) = \tfrac12 \sum_{i=1}^{4n} c^N(IA,e_i,I e_i),
\end{equation*}
where \( \{e_i\} \) is an orthonormal basis for \( TN \).  This one-form is
independent of~\( I \) and globally defined~\cite{Ivanov:QKT}.

\begin{proposition}
  Let \( N \)~be a \QKT manifold.  If \( \dim N = 4 \), suppose also that
  \( N \)~is self-dual.  Then \( N \) is of instanton type if and only if
  its torsion one-form~\( \tau \) satisfies \(
  d\tau\in\bigcap_{I}\Lambda^{1,1}_I \).
\end{proposition}

\begin{proof}
  As \( \nabla \) is a metric connection, we have that \( \Nq \vol^N =
  \xi\cdot\vol^N = \sum_{i=1}^{4n} g^N(\xi e_i,e_i) \vol^N \).  Let \(
  \{f_1,\dots,f_n\} \) be an orthonormal quaternionic basis for~\( TN \).
  Then
  \begin{equation*}
    \begin{split}
      &\sum_{i=1}^{4n} g^N(\xi_A e_i,e_i)
      = 4 \sum_{j=1}^n g^N(\xi_A f_j,f_j)\\
      &= 4 \sum_{j=1}^n \tfrac16
      \bigl(
      -c^N(IA,f_j,If_j)-c^N(JA,f_j,Jf_j)-c^N(KA,f_j,Kf_j)
      \bigr)\\
      &\qquad - \tfrac1{12}
      \bigl(
      c^N(JA,Kf_j,If_j)+c^N(KA,If_j,Jf_j)+c^N(IA,Jf_j,Kf_j)\\
      &\qquad\qquad
      -c^N(KA,Jf_j,If_j)-c^N(IA,Kf_j,Jf_j)-c^N(JA,If_j,Kf_j)
      \bigr)\\
      &=-\tfrac16\sum_{j=1}^n
      c^N(IA,f_j,If_j)+c^N(IA,Jf_j,Kf_j)+c^N(JA,f_j,Jf_j)\\
      &\qquad\qquad +c^N(JA,Kf_j,If_j)+c^N(KA,f_j,Kf_j)+c^N(KA,If_j,Jf_j)\\
      &=-\tfrac14\tau(A).
    \end{split}
  \end{equation*}
  Thus the curvature of \( \kappa^{\mathbb R} \) is \( -\tfrac14 d\tau \),
  giving the result.
\end{proof}

\section{An Inverse Construction}
\label{sec:inverse}
Let \( N^{4n} \) be a \QKT manifold of instanton type.  We follow the
constructions and notation of \cite{Swann:MathAnn}.  Set \( \UM(N) \) to be
the bundle
\begin{equation*}
  \UM(N) = P \times_{\SP(n)\SP(1)} \mathbb H^*/\{\pm 1\},
\end{equation*}
where \( P \) is the principal \( \SP(n)\SP(1) \)-bundle of frames over~\(
N \).  This is a real line bundle over the bundle \( \mathcal S(N) \) of
quaternion bases for~\( \mathcal G \).  Let \( \omega \)~be the connection
one-form on~\( P \) defined by the torsion-connection~\( \nabla \).  As \(
\nabla \)~is an \( \SP(n)\SP(1) \)-connection, \( \omega \)~takes values
in~\( \sP(n)\oplus\sP(1) \).  Let \( \omega_- \) be the \(
\sP(1)\cong\im\mathbb H \) part of~\( \omega \).  The form \( \omega_-
\)~defines a splitting of \( T\UM(N) = \mathcal H + \mathcal V \), with \(
\mathcal H\cong T_pN \) and \( \mathcal V\cong\mathbb H \).  The projection
to~\( \mathcal S(N) \) defines three almost complex structures on~\(
\mathcal H \).  Combining these with those on~\( \mathcal V \) one gets a
hypercomplex structure on~\( \UM(N) \).  This follows from the results
of~\cite{Pedersen-PS:hypercomplex}, since Proposition~\ref{prop:Nq} shows
that \( \omega_- \) extends to a torsion-free quaternionic connection and
\( \UM(N) \) is \( \UM_{(n-1)/n}(N) \) twisted by the instanton bundle \(
(\kappa^{\mathbb R})^{-1/4(n+1)} \).

Consider the function 
\begin{equation*}
  \mu=x\bar x 
\end{equation*}
on~\( \UM(N) \), where \( x \)~is the quaternionic coordinate on \( \mathbb
H^* \).  Then \( d\mu = x\bar \psi + \psi \bar x \), where \( \psi =
dx - x\omega_- \), and so \( I d\mu = -x\bar \psi i +i \psi \bar x \).
Thus we have
\begin{equation*}
  dId\mu = -\psi\wedge\bar \psi i - i\psi\wedge\bar \psi -
  x\Omega_-\bar xi -i x\Omega_-\bar x,
\end{equation*}
where \( \Omega_- = d\omega_-+\omega_-\wedge\omega_- = -\tfrac12 (\beta_Ii
+ \beta_Jj + \beta_Kk) \) is the \( \sP(1) \)-curvature.

The function \( \mu \) is an \HKT potential if \( F_I=\tfrac12(dId-JdId)\mu
\) defines a Kähler form for~\( I \) associated to a non-degenerate metric
via \( g(I\cdot,\cdot) = F_I(\cdot,\cdot) \).  Vertically, \( F_I \)~is
given by the \( \psi \) terms and is the Kähler form of the flat structure
on~\( \mathbb H^* \).

Horizontally, we get a condition on~\( \Omega_- \).  At \( x=1 \) the
horizontal part of~\( F_I \) is \( \tfrac12(\beta_I - J\beta_I) \).  As \(
N \) is of instanton type we have \( \beta_I=\sigma_I^{1,1}(\cdot,I\cdot )
\), in the notation of the previous section.  This gives that the
horizontal part of~\( g \) at~\( x=1 \) is \(
-\tfrac12(1+J_N)\sigma_I^{1,1} = \sigma^q \), where \( \sigma^q \) is the
component of the symmetric part of~\( -2B \) that is of type \( (1,1) \)
for each \( I_N \).  For general \( x \), the horizontal part of~\( g \) is
now \( x\bar x\sigma^q \).  We call \( \sigma^q =
\tfrac12(\beta_I-J_N\beta_I)(\cdot,I_N\cdot) \) the \emph{curvature metric}
of~\( N \).  Recall that \( \beta_I \)~is defined by
equation~\eqref{eq:R-I} and so we have
\begin{equation*}
  \sigma^q = \frac1{4n}\sum_{i=1}^{4n}
  R^N(X,I_NY,e_i,I_Ne_i)+R^N(J_NX,K_NY,e_i,I_Ne_i). 
\end{equation*}

The group~\( \mathbb H^* \) acts on~\( \UM(N) \) via left-multiplication
and generates a special homothety of type~\( (2,-2) \), i.e., we have a
\Dsym[-2]-symmetry.  Using the results of~\S\ref{sec:geography} we now
obtain \HKT structures with \Dsym-symmetry for other values of~\( \alpha
\).

\begin{theorem}
  \label{thm:UM}
  Let \( N \) be a \QKT manifold of instanton type whose curvature metric
  is non-degenerate of signature~\( (4p,4q) \).  Then \( \UM(N) \) carries
  an \HKT structures with \Dsym-symmetry for each~\( \alpha\ne0 \).  For \(
  \alpha<0 \) the metric has signature~\( (4p+4,4q) \); for \( \alpha>0 \)
  the signature is~\( (4q+4,4p) \). \qed
\end{theorem}

\begin{example}
  If \( N \) is a quaternionic Kähler manifold with non-zero scalar
  curvature~\( s \), then the curvature metric is \( \sigma^q = s' g^N \),
  where \( s' \)~is a positive dimension-dependent multiple of~\( s \).
  The \HKT structure on \( \UM(N) \) is thus positive definite exactly when
  \( s\alpha <0 \).  In \cite{Swann:MathAnn} it was shown that the metric
  constructed from \( \mu=x\bar x \) is hyperKähler with
  \Dsym[-2]-symmetry.  For \( N=\HP(n) \), we have \( \UM(N)=\mathbb
  H^{n+1}/\{\pm1\} \) with the flat hyperKähler metric.  For \( N=\HH(n)
  \), or any other quaternionic Kähler manifold of negative scalar
  curvature, the hyperKähler metric on~\( \UM(N) \) is indefinite; but
  Theorem~\ref{thm:UM} above shows that \( \UM(N) \) has a definite \HKT
  metric with \Dsym[1]-symmetry.
\end{example}

\begin{remark}
  \label{rem:descent}
  Suppose \( N \) is a \QKT manifold constructed as the \( \mathbb H^*
  \)-quotient of an \HKT manifold with \Dsym-symmetry as in
  Theorem~\ref{thm:H-quotient}.  The proof of Proposition~\ref{prop:sp1}
  shows that \( \sigma^q \) is determined by the \HKT potential and is the
  horizontal part of~\( g \) up to a constant.  But that is also how \( g^N
  \) is constructed.  So \( N \) has curvature metric~\( \sigma^q \)
  proportional to~\( g^N \) and the construction in this section is inverse
  to that of Theorem~\ref{thm:H-quotient}.
\end{remark}

\begin{remark}
  In \cite{Howe-OP:QKT} a construction similar to Theorem~\ref{thm:UM} is
  given, with different assumptions on the base and with the conclusion
  that \( \UM(N) \) is hyperKähler.  The essential condition
  in~\cite{Howe-OP:QKT} is that \( dc^N \) should be of type~\( (2,2) \)
  for each~\( I \).  The appendix of that paper contains a proof that this
  implies that the curvature metric is quaternionic Kähler.
\end{remark}

Ivanov~\cite{Ivanov:QKT} showed that every metric \( \bar g^N = e^u g^N \)
conformal to a \QKT metric~\( g^N \) admits a \QKT connection \( \bar
\nabla^N \).  The torsion-one form satisfies \( \bar\tau = \tau - (2n+1)du
\), so \( (\bar g^N,\bar \nabla^N) \) is of instanton type if and only if
\( (g^N,\nabla^N) \) is.  One computes that 
\begin{equation*}
  \begin{split}
    \bar\beta_I(X,Y) = \beta_I(X,Y) &+ (\Nq du)(IX,Y) - (\Nq du)(IY,X)\\
    &- (Jdu\wedge Kdu) (X,Y)
  \end{split}
\end{equation*}
and hence
\begin{equation}
  \label{eq:conformal-change}
  \bar \sigma^q = \sigma^q + \tfrac12 (1+I+J+K)\nabla^qdu + \tfrac12
  (d^{\mathbb H}u)^2.
\end{equation}
We may use this result in several ways.  Firstly, note that it shows that
for a general \QKT manifold~\( N \) of instanton type, the curvature
metric~\( \sigma^q \) need not be proportional to~\( g^N \).  Secondly, if \(
g^N \) has degenerate curvature metric, then in a neighbourhood of any
point we can choose~\( u \) so that \( \bar \sigma^q \)
in~\eqref{eq:conformal-change} is non-degenerate.

\begin{proposition}
  Suppose \( N \) is a quaternionic manifold.  Then locally, \( N \)~admits
  a positive definite \QKT structure.
\end{proposition}

\begin{proof}
  Fix a volume form~\( \vol_0 \) on~\( N \).  Then there is a unique
  torsion-free quaternionic connection~\( \nabla^0 \) on~\( N \) such that
  \( \nabla^0\vol_0=0 \).  If the curvature metric \( \sigma_0 \) is not
  positive definite, replace \( \vol_0 \) be \( e^{2nu}\vol_0 \) for some
  function~\( u \); then \( \sigma_0 \) will change as
  in~\eqref{eq:conformal-change} and we may choose~\( u \) so that the
  curvature metric is positive definite in a neighbourhood of a given
  point.  Now set \( g^N=\sigma^0 \).  Since \( \nabla^0 \) preserves~\(
  \vol_0 \) the bundle \( \kappa^{\mathbb R} \) is an instanton and we may
  use \( g^N \) and \( \nabla^0 \) to construct an \HKT structure on~\(
  \UM(N) \) as in Theorem~\ref{thm:UM}.  By Theorem~\ref{thm:H-quotient} we
  get a \QKT structure on~\( N \) which by Remark~\ref{rem:descent} has \(
  g^N \) as its metric, up to a constant scale.
\end{proof}

\providecommand{\bysame}{\leavevmode\hbox to3em{\hrulefill}\thinspace}

\end{document}